\RequirePackage{fix-cm}
\documentclass[smallextended]{svjour3}       
\smartqed  
\usepackage{graphicx}
\usepackage{amsmath}
\usepackage{amssymb,latexsym}
\usepackage{oldlfont}
\usepackage{amsfonts}
\begin{document}

\title{Saddle representations of positively homogeneous functions by linear functions
}

\titlerunning{Saddle representations of p.h. functions}        

\author{Valentin V. Gorokhovik        \and
        Marina Trafimovich}


\institute{Valentin V. Gorokhovik,  Corresponding author  \at
             Institute of Mathematics,
             National Academy of Sciences of Belarus,\\
             ul. Surganova, 11, 220072 Minsk, Belarus\\
             Tel.: +375-17-2841951\\
              Fax: +375-17-2841701\\
            \email{gorokh@im.bas-net.by}
           \and
          Marina Trafimovich  \at
              Sukhoi State Technical University of Gomel,\\
              Pr. Octiabria, 48, 246746, Gomel, Belarus\\
              \email{marvoitova@tut.by}
}

\date{Received: date / Accepted: date}

\maketitle

\begin{abstract}
We say that a positively homogeneous function admits a saddle representation by linear functions iff it admits both an inf-sup-representation and a sup-inf-representation with the same two-index family of linear functions.
In the paper we show that each continuous positively homogeneous function
can be associated with a two-index family of linear functions which provides its saddle representation.
We also establish characteristic properties of those two-index families of linear functions which provides saddle representations of functions belonging to
the subspace of Lipschitz continuous positively homogeneous functions as well as the subspaces of difference sublinear and piecewise linear functions.
\keywords{Positively homogeneous functions  \and Saddle representation \and Lipschitz continuity \and Difference sublinear functions \and Piecewise linear functions}

\subclass{49J52 \and 54C35}
\end{abstract}

\section{Introduction}
\label{intro}
Positively homogeneous (p.h.) functions play a crucial role in nonsmooth analysis and optimization. Let us recall that directional derivatives widely used as first order approximation of nonsmooth functions in different theories of generalized differentiation (see, for example, Demyanov and Rubinov \cite{DemRub90,DemRub95}, Mordukhovich \cite{Mord1,Mord2}, Pallaschke and Rolewich \cite{PalRol}, Rockafellar and Wets \cite{RW98} and references therein) are p.h. functions. Besides, nonlinear programming problems with p.h. objective and constraint functions are used as approximations of initial nonsmooth nonlinear programming problems in deriving optimality conditions \cite{BonSha} and in numerical methods \cite{Bagirov}. Moreover, optimization problems with homogeneous data are of independent interest as an important  class of mathematical models describing real problems (see, for instance, \cite{Rub} and \cite{HU}). Both necessary and sufficient optimality conditions in nonsmooth optimization problems are often formulated as a condition of nonnegativity or positivity of a some p.h. function formed by directional derivatives of functions involved in the optimization problem under consideration. At the same time, since a p.h. function is uniquely defined by its restriction to the unit sphere, it may be so much complicated as an arbitrary function defined on the unit sphere. It shows that the verification of properties (even such simple ones as nonnegativity or positivity) of p.h. function may be rather difficult in general. To overcome these difficulties we need   to study the structure of most important classes of p.h. functions, in particular, the possibilities of their representation through linear functions.

The paper deals with the representations of p.h. functions by families of linear functions using consecutive operations of the pointwise infimum and the pointwise supremum. Such representations were studied by Demyanov \cite{Dem99,DR2000}, Castellani \cite{Cast1,Cast2}, Castellani and Uderzo \cite{Cast3}, Gorokhovik and Starovoitova \cite{GS2011}. The novelty of the results presented in this note is that we are interested in existing such two-index families of linear functions which provides both an inf-sup-representation and a sup-inf-representation of a p.h. function simultaneously; such representations are referred below as saddle ones. The advantage of  saddle representations in comparison with  inf-sup-representations and sup-inf-representations is their universality: for example, the same saddle representation of a p.h. function can be used for to derive dual characteristics of both its nonnegativity and its nonpositivity as well as for to find both a direction of steepest descent and a direction of steepest ascent.

\section{Preliminaries}
\label{sec:1}

Recall that a function $p: {\mathbb R}^n \to {\mathbb R}$ is said to be \textit{positively homogeneous} (of first degree) iff $p(\lambda x) = \lambda p(x)$ for all $x \in {\mathbb R}^n$ and all positive reals $\lambda > 0.$

\vspace{2mm}

With respect to standard algebraic operations, the collection of p.h. functions, defined on ${\mathbb R}^n,$ forms a real vector space denoted below by ${\mathcal P}({\mathbb R}^n).$ In addition, both the pointwise maximum and the pointwise minimum of a finite subfamily of p.h. functions are also p.h. functions. These order operations are consistent with algebraic operations in such a way that ${\mathcal P}({\mathbb R}^n)$ is a vector lattice.

Here, we restrict ourselves with consideration of p.h. functions belonging to a number of subspaces of ${\mathcal P}({\mathbb R}^n),$ in particular, to the subspace ${\mathcal P}_C({\mathbb R}^n)$ of continuous p.h. func\-tions, the subspace ${\mathcal P}_L({\mathbb R}^n)$ of Lipschitz continuous p.h. func\-tions, the subspace ${\mathcal P}_{DC}({\mathbb R}^n)$ of difference sublinear functions, the subspace ${\mathcal PL}({\mathbb R}^n)$ of piecewise linear functions, and the subspace ${\mathcal L}({\mathbb R}^n)$ of linear functions. These subspaces are connected to each other with the following chain of inclusions
\begin{equation}\label{e1}
{\mathcal L}({\mathbb R}^n) \subset {\mathcal PL}({\mathbb R}^n) \subset {\mathcal P}_{DC}({\mathbb R}^n) \subset {\mathcal P}_{L}({\mathbb R}^n) \subset {\mathcal P}_{C}({\mathbb R}^n) \subset {\mathcal P}({\mathbb R}^n).
\end{equation}
Each of the above subspaces, except the subspace of linear functions ${\mathcal L}({\mathbb R}^n),$ contains both the pointwise maximum and the pointwise minimum of a finite subfamily of functions belonging to it and, consequently, is a vector sublattice of ${\mathcal P}({\mathbb R}^n).$

Apart from the aforementioned subspaces of p.h. functions, we will consider the classes of sublinear (convex and positively homogeneous) and superlinear (concave and positively homogeneous) functions denoted below by ${\mathcal P}_{subl}({\mathbb R}^n)$ and  ${\mathcal P}_{supl}({\mathbb R}^n)$ respectively. These classes are convex cones in ${\mathcal P}_{DC}({\mathbb R}^n).$ Moreover, both ${\mathcal P}_{subl}({\mathbb R}^n)$ and  ${\mathcal P}_{supl}({\mathbb R}^n)$ are generating for ${\mathcal P}_{DC}({\mathbb R}^n),$ that is \linebreak ${\mathcal P}_{DC}({\mathbb R}^n)= {\mathcal P}_{subl}({\mathbb R}^n) - {\mathcal P}_{subl}({\mathbb R}^n)$ and ${\mathcal P}_{DC}({\mathbb R}^n)= {\mathcal P}_{supl}({\mathbb R}^n) - {\mathcal P}_{supl}({\mathbb R}^n).$ In addition,
$
{\mathcal P}_{subl}({\mathbb R}^n) = - {\mathcal P}_{supl}({\mathbb R}^n)
\,\,\text{and}\,\,
{\mathcal P}_{subl}({\mathbb R}^n) \cap {\mathcal P}_{supl}({\mathbb R}^n) = {\mathcal L}({\mathbb R}^n).
$

Due to the isomorphism between sublinear functions and compact convex sets (see, for instance, \cite{HULem1}) each sublinear function $\varphi:{\mathbb R}^n \to {\mathbb R}$ can be uniquely associated with the convex compact subset
$$
\partial \varphi := \{a \in {\mathbb R}^n\,\,: \,\,\langle a, x \rangle \le \varphi(x)\,\,\forall\,\,x \in {\mathbb R}^n\},
$$
such that
$$
\varphi(x) = \max\limits_{a \in \partial \varphi}\langle a, x \rangle\,\,\forall\,\,x \in {\mathbb R}^n.
$$
Thus, each sublinear function $\varphi$ is the pointwise maximum (upper envelope) of the collection of its linear minorants. The set $\partial \varphi$ is called the \textit{(lower) subdifferential} of $\varphi$ (at the null point).

Likewise, each superlinear function $\psi:{\mathbb R}^n \to {\mathbb R}$ can be uniquely associated with the convex compact subset
$$
\partial^+ \psi := \{a \in {\mathbb R}^n\,\,: \,\,\langle a, x \rangle \ge \psi(x)\,\,\forall\,\,x \in {\mathbb R}^n\},
$$
called the upper subdifferential of $\psi$ (at the null point), such that
$$
\psi(x) = \min\limits_{a \in \partial^+ \psi}\langle a, x \rangle\,\,\forall\,\,x \in {\mathbb R}^n\}.
$$
Consequently, each superlinear function $\varphi$ is the pointwise minimum (lower envelope) of the collection of its linear majorants.

Now, let us consider a family of sublinear functions $\Phi \subset {\mathcal P}_{subl}({\mathbb R}^n)$ and assume that
$$
\inf\limits_{\varphi \in \Phi}\varphi(x) > -\infty\,\,\forall\,\,x \in {\mathbb R}^n.
$$
Since sublinear functions are continuous on ${\mathbb R}^n$, we have due to \cite[p.~300]{Bourb} that the function
$$
p(x):= \inf\limits_{\varphi \in \Phi}\varphi(x)\,\,\forall\,\,x \in {\mathbb R}^n
$$
is upper semicontinuous and positively homogeneous on ${\mathbb R}^n.$

Likewise, whenever a family of superlinear functions $\Psi \subset {\mathcal P}_{supl}({\mathbb R}^n)$  holds
$$
\sup\limits_{\psi \in \Psi}\varphi(x) < +\infty\,\,\forall\,\,x \in {\mathbb R}^n,
$$
the function
$$
q(x) := \sup\limits_{\psi \in \Psi}\varphi(x)\,\,\forall\,\,x \in {\mathbb R}^n
$$
is lower semicontinuous and positively homogeneous on ${\mathbb R}^n.$

Demyanov and Rubinov \cite{DemRub90,DemRub95,Dem82} proved the converse statement: any upper (lower) semicontinuous p.h. function $p:{\mathbb R}^n \to {\mathbb R}$ is the pointwise infimum (supremum) of some family of sublinear (superlinear) functions. These results were extended firstly by Uderzo \cite{Ud2000} to infinite-dimensional uniformly convex Banach spaces and then by Gorokhovik \cite{Gor17} to arbitrary normed spaces. 

\begin{proposition}{\rm \cite{Dem82,Ud2000,Gor17}}\label{prop}
Let $p:{\mathbb R}^n \to {\mathbb R}$ be a positively homogeneous function. Then
\begin{itemize}
\item[\hspace*{2pt}$(i)$]
 $p: {\mathbb R}^n \to {\mathbb R}$ is upper semicontinuous on ${\mathbb R}^n$ if and only if there exists a family $\Phi$ of sublinear functions such that
\begin{equation}\label{e2}
  p(x) = \inf\limits_{\varphi \in \Phi}\max\limits_{a \in \partial \varphi}\langle a, x \rangle\,\,\forall\,\,x \in {\mathbb R}^n;
\end{equation}
\vspace{2pt}
\item[$(ii)$]
$p: {\mathbb R}^n \to {\mathbb R}$ is lower semicontinuous on ${\mathbb R}^n$ if and only if there exists a family $\Psi$ of superlinear functions such that
\begin{equation}\label{e3}
  p(x) = \sup\limits_{\psi \in \Psi}\min\limits_{a \in \partial^+ \psi}\langle a, x \rangle\,\,\forall\,\,x \in {\mathbb R}^n;
\end{equation}
\vspace{2pt}
\item[$(iii)$]
$p: {\mathbb R}^n \to {\mathbb R}$ is continuous on ${\mathbb R}^n$ if and only if there exist both a family $\Phi$ of sublinear functions and a family $\Psi$ of superlinear functions which hold the equality \eqref{e2} and the equality \eqref{e3}, respectively.
\end{itemize}
\end{proposition}

It follows from representations \eqref{e2} and \eqref{e3} that lower and upper semicontinuous p.h. functions can be ``constructed'' from families of linear functions by means of consecutive operations of the pointwise infimum and the pointwise supremum.
For continuous p.h. functions the both representations, \eqref{e2} and \eqref{e3}, hold.

We say that a p.h function $p:{\mathbb R}^n \to {\mathbb R}$ admits \textit{an inf-sup-representation} by linear functions iff there exists a two-index family $\{b_{is} \in {\mathbb R}^n\,|\,i \in I,\,s \in S\}$  (here and throughout below $I$ and $S$ are index sets) such that
\begin{equation}\label{e4}
p(x) = \inf\limits_{i \in I}\sup\limits_{s \in S}\langle b_{is},\, x \rangle \,\,\forall\,\,x \in {\mathbb R}^n.
\end{equation}
Likewise, we say that a p.h function $p:{\mathbb R}^n \to {\mathbb R}$ admits \textit{a sup-inf-representation} by linear functions iff there exists a two-index family $\{c_{is} \in {\mathbb R}^n\,|\,i \in I,\,s \in S\}$ such that
\begin{equation}\label{e5}
p(x) = \sup\limits_{s \in S}\inf\limits_{i \in I}\langle c_{is},\, x \rangle
\,\,\forall\,\,x \in {\mathbb R}^n.
\end{equation}

It follows from Proposition \ref{prop} that each upper (lower) semicontinuous p.h. function admits an inf-sup-representation (a sup-inf-representation) by linear functions while each continuous p.h. function admits  both an inf-sup-representation and a sup-inf-representation by linear functions. It is worth observing that due to Proposition \ref{prop} an inf-sup-representation and a sup-inf-representation of the same continuous p.h. function $p$ are generally provided with different two-index families $\{b_{is} \in {\mathbb R}^n\,|\,i \in I,\,s \in S\}$ and $\{c_{is} \in {\mathbb R}^n\,|\,i \in I,\,s \in S\}$.

In some cases, it is more convenient to use an inf-sup-representation, while in other cases a sup-inf-representation is more convenient. For example, the usage of an inf-sup-representation of continuous positively homogeneous function is more convenient than its sup-inf-representation when we need to characterize nonnegativity of this function, and conversely, for characterization of nonpositivity a sup-inf-representation is more convenient. Motivated by these circumstances Demyanov \cite{DR2000} for Lipshitz continuous positively homogeneous functions and Gorokhovik and Trafimovich \cite{GT2015} for continuous ones developed methods for converting a inf-sup representation into a sup-inf-representation and vice-versa.

We say that a p.h function $p:{\mathbb R}^n \to {\mathbb R}$ admits \textit{a saddle representation} by linear functions iff there exists a two-index family $\{a_{is} \in {\mathbb R}^n\,|\,i \in I,\,s \in S\}$ such that
\begin{equation}\label{e6}
p(x) = \inf\limits_{i \in I}\sup\limits_{s \in S}\langle a_{is},\, x \rangle = \sup\limits_{s \in S}\inf\limits_{i \in I}\langle a_{is},\, x \rangle, \,\,\forall\,\,x \in {\mathbb R}^n,
\end{equation}
that is, iff there exists a two-index family $\{a_{is} \in {\mathbb R}^n\,|\,i \in I,\,s \in S\}$ which provides both {an inf-sup-representation} and {a sup-inf-representation} of $p$ by linear functions.

The main result of the paper is the proof that we can associate with each continuous p.h. function $p$ such two-index family $\{a_{is} \in {\mathbb R}^n\,|\,i \in I,\,s \in S\}$ which provides a saddle representation of $p$ by linear functions. We also establish additional characteristic properties of those two-index families which provides saddle representations of Lipschitz continuous p.h. functions as well as difference-sublinear and piecewise-linear ones.

\section{Main results}

\begin{theorem}\label{th1}
For a function $p:{\mathbb R}^n \to {\mathbb R}$ to be positively homogeneous and continuous on ${\mathbb R}^n$, it is necessary and sufficient that there exists a two-index family $\{a_{is} \in {\mathbb R}^n\,|\,i \in I,\,s \in S\}$ which satisfies, for each $x \in {\mathbb R}^n$, the conditions
\begin{equation}\label{e7}
-\infty < \inf\limits_{i \in I}\langle a_{is},\, x \rangle\,\,\,\forall\,\,s \in S\,\,\,\,
\text{and}\,\,\,\,\sup\limits_{s \in S}\langle a_{is},\, x \rangle < +\infty\,\,\,\forall\,\,i \in I,
\end{equation}
and provides the saddle representation \eqref{e6} of $p.$
\end{theorem}

{\it Proof}  Let $p:{\mathbb R}^n \to {\mathbb R}$ be continuous p.h. function and let $\Phi=\{\varphi_i:i\in I\}$ and $\Psi = \{\psi_s:s \in S\}$ be families of sublinear and superlinear functions, which satisfy equalities \eqref{e2} and \eqref{e3}, respectively. The existence  of such families is provided by the statement $(iii)$ of Proposition \ref{prop}. For any given $i \in I$ and $s \in S$ we have from the equalities \eqref{e2} and \eqref{e3} that
$$
\psi_s(x) \le p(x) \le \varphi_i(x)\,\,\forall\,\,x \in {\mathbb R}^n.
$$
Since the function $\psi_s$ is superlinear and the function $\varphi$ is sublinear, it follows from separation theorems of convex sets \cite{Rock} that there exists a linear function $x \to \langle a_{is},\,x \rangle,$ where $a_{is} \in {\mathbb R}^n,$ such that
\begin{equation}\label{e8}
    \psi_s(x) \le \langle a_{is},\,x \rangle \le \varphi_i(x)\,\,\forall\,\,x \in {\mathbb R}^n.
\end{equation}
Taking the infimum over $i \in I$ in the last inequalities, we obtain
\begin{equation}\label{e9}
    \psi_s(x) \le \inf\limits_{i \in I}\langle a_{is},\,x \rangle \le \inf\limits_{i \in I}\varphi_i(x)\,\,\forall\,\,x \in {\mathbb R}^n.
\end{equation}
It shows that
$$
-\infty < \inf\limits_{i \in I}\langle a_{is},\,x \rangle\,\,\forall\,\,s \in S.
$$
Now, taking the supremum over $s \in S$ in the inequalities \eqref{e9}, we get
 $$
p(x) = \sup\limits_{s \in S}\psi_s(x) \le \sup\limits_{s \in S}\inf\limits_{i \in I}\langle a_{is},
\, x \rangle \le \inf\limits_{i \in I}\varphi_i(x) = p(x)\,\,\forall\,\,x \in {\mathbb R}^n
$$
Consequently,
$$
p(x) = \sup\limits_{s \in S}\inf\limits_{i \in I}\langle a_{is},\, x \rangle \,\,\forall\,\,x \in {\mathbb R}^n.
$$

If we take firstly the supremum over $s \in S$ in the inequalities \eqref{e8} and then the infimum  over $i \in I$, we prove that
$$
\sup\limits_{s \in S}\langle a_{is},\,x \rangle < +\infty\,\,\forall\,\,i \in I,
$$
and
$$
p(x) = \inf\limits_{i \in I}\sup\limits_{s \in S}\langle a_{is},\, x \rangle \,\,\forall\,\,x \in {\mathbb R}^n.
$$
The necessary  part of Theorem \ref{th1} is proved.

To prove the sufficient part we first observe that a function $p:{\mathbb R}^n \to {\mathbb R}$ which admits an inf-sup-representation or a sup-inf-representation by linear functions is positively homogeneous.

Besides, an arbitrary function $p:{\mathbb R}^n \to {\mathbb R}$ admitting  an inf-sup-representa\-tion $p(x) = \inf\limits_{i \in I}\sup\limits_{s \in S}\langle a_{is},\, x \rangle \,\,\forall\,\,x \in {\mathbb R}^n$ with such two-index family $\{a_{is} \in {\mathbb R}^n\,|\,i \in I,\,s \in S\}$ which satisfies, for each $x \in {\mathbb R}^n$, the condition
$\sup\limits_{s \in S}\langle a_{is},\, x \rangle < +\infty\,\,\,\forall\,\,i \in I,$ is upper semicontinuous on ${\mathbb R}^n.$ Really, in this case $p$ is the lower envelope of the family of sublinear functions $\{\varphi_i: x \to \sup\limits_{s \in S}\langle a_{is},\, x \rangle,\,i \in I\}$ each of which takes finite values for all $x \in {\mathbb R}^n.$ Because real-valued sublinear functions $\varphi_i,\,i\in I,$ are continuous on ${\mathbb R}^n,$ their lower envelope, that is the function $p,$ is upper semicontinuous on ${\mathbb R}^n.$

Likewise, a function $p:{\mathbb R}^n \to {\mathbb R}$ admitting  an sup-inf-representation $p(x) = \sup\limits_{s \in S}\inf\limits_{i \in I}\langle a_{is},\, x \rangle$  $\forall\,\,x \in {\mathbb R}^n$ with such two-index family $\{a_{is} \in {\mathbb R}^n\,|\,i \in I,\,s \in S\}$ which satisfies, for each $x \in {\mathbb R}^n$, the condition
$- \infty < \inf\limits_{i \in I}\langle a_{is},\, x \rangle\,\,\,\forall\,\,i \in I,$ is lower semicontinuous on ${\mathbb R}^n.$
\qed

\begin{remark}\label{r1}
It follows from the proof of Theorem \ref{th1} that the two-index family \linebreak $\{a_{is} \in {\mathbb R}^n\,|\,i \in I,\,s \in S\}$ which provides the saddle representation for the given continuous positively homogeneous function $p:{\mathbb R}^n \to {\mathbb R}$ is nonuniquely defined. The nonuniqueness arises from the nonunique choice of exhaustive families of upper convex and lower concave approximations of the function $p$ and, in addition, from the possible nonuniqueness of a linear function $x \to \langle a_{is},\,x \rangle$ which separates an upper convex approximation $\varphi_i$ and a lower concave one $\psi_s$ in \eqref{e8}.
\end{remark}

\begin{remark}\label{r2}
Let $\{a_{is} \in {\mathbb R}^n\,|\,i \in I,\,s \in S\}$ be the two-index family of linear functions providing a saddle representation for the p.h.  function $p$ that was constructed in the proof of Theorem 1. It is evident that the cardinalities of the index sets $I$ and $S$ are equal to those of the exhaustive families of upper convex and lower concave approximations which were used for constructing this family. Consequently, the less the cardinalities of the exhaustive families, the less those of $I$ and $S$. Clearly, we are interested in reducing the cardinalities. To this end we can use the methods of reducing exhaustive families which were developed in \cite{Roshch,GPU2010,GKKU15}.
\end{remark}

For the next theorem we need the following characterizations of Lipschitz continuous positively homogeneous functions.

\begin{proposition}{\rm \cite{GS2011}}\label{pr2}
The following statements are equivalent
\begin{itemize}
\item[\hspace*{2pt}$(i)$]
$p: {\mathbb R}^n \to {\mathbb R}$ is positively homogeneous and Lipschitz continuous on ${\mathbb R}^n;$
\vspace{2pt}
\item[$(ii)$]
there exists a uniformly bounded family $\Phi$ of sublinear functions which satisfies the equality \eqref{e2};
\vspace{2pt}
\item[$(iii)$]
there exists a uniformly bounded family $\Psi$ of superlinear functions which satisfies the equality \eqref{e3}.
\end{itemize}
\end{proposition}

The family $\Gamma$ of positively homogeneous functions is called \textit{uniformly bounded} if there exists a constant $M > 0$ such that $|p(x)| \le M\|x\|$ for all $x \in {\mathbb R}^n$ and all $p \in \Gamma.$

We also note that for a sublinear (respectively, superlinear) function $p$ the condition $|p(x)| \le M\|x\|$ $\forall\,\,x \in {\mathbb R}^n$ is equivalent to $p(x) \le M\|x\|\,\,\forall\,\,x \in {\mathbb R}^n$ (respectively, to $-M\|x\| \le p(x)\,\,\forall\,\,x \in {\mathbb R}^n$).

\begin{theorem}\label{th2}
For a function $p:{\mathbb R}^n \to {\mathbb R}$ to be positively homogeneous and Lipschitz continuous on ${\mathbb R}^n$, it is necessary and sufficient that there exists a two-index family $\{a_{is} \in {\mathbb R}^n\,|\,i \in I,\,s \in S\}$ which satisfies, for some fixed $M > 0$ and each $x \in {\mathbb R}^n$, the conditions
\begin{equation}\label{e10}
-M\|x\| < \inf\limits_{i \in I}\langle a_{is},\, x \rangle\,\,\,\forall\,\,s \in S\,\,\,\,
\text{and}\,\,\,\,\sup\limits_{s \in S}\langle a_{is},\, x \rangle < M\|x\|\,\,\,\forall\,\,i \in I,
\end{equation}
and provides the saddle representation \eqref{e6} of $p$.
\end{theorem}

{\it Proof} Let $p:{\mathbb R}^n \to {\mathbb R}$ be a Lipshitz continuous positively homogeneous function. Then, due to Proposition~\ref{pr2}, we can choose a uniformly bounded family $\Phi$ of sublinear functions as well as  a uniformly bounded family $\Psi$ of superlinear functions which satisfy the equalities \eqref{e2} and \eqref{e3}, respectively. Repeating the arguments of the proof of Theorem \ref{th2} with the chosen families $\Phi$ and $\Psi$ we construct a two-index family $\{a_{is},\,i \in I,\,s \in S\},$ which satisfies the conditions\eqref{e10} and provides for $p$ the saddle representation \eqref{e6}.

Conversely, it follows from \eqref{e6} and \eqref{e10} that the family of superlinear functions $\psi_s:x \to \inf\limits_{i \in I}\langle a_{is},\,x\rangle,\,s \in S,$ is uniformly bounded and satisfies the equality \eqref{e3}.  Hence,  the function $p$ is positively homogeneous and, due to Proposition \ref{pr2}, it is Lipschitz continuous.
\qed

\begin{theorem}\label{th3}
A function $p:{\mathbb R}^n \to {\mathbb R}$ is difference sublinear if and only if there exist two bounded families of vectors in ${\mathbb R}^n$, $\{b_{s} \in {\mathbb R}^n\,|\,s \in S\}$ and $\{c_{i} \in {\mathbb R}^n\,|\,i \in I\},$ such that the two-index family $\{a_{is} :=b_s - c_i \in {\mathbb R}^n\,|\,i \in I,\,s \in S\}$ generated by this families provides the saddle representation \eqref{e6} of $p$ by linear functions.
\end{theorem}

{\it Proof} A function $p:{\mathbb R}^n \to {\mathbb R}$ is difference sublinear if and only if $p$ can be represented as a difference of two sublinear functions. Let $p(x) = \varphi_1(x) - \varphi_2(x)\,\,\forall\,\, x \in {\mathbb R}^n,$
where $\varphi_1,\,\varphi_2 \in {\mathcal P}_{subl}({\mathbb R}^n),$ be some such representation of $p.$ Consider the subdifferentials $\partial \varphi_1$ and $\partial \varphi_2$ of the functions $\varphi_1$ and $\varphi_2$, respectively, and denote by $\{b_s\,|\,s \in S\}$ the set of exposed points of the subdifferential $\partial \varphi_1$ indexed by elements of $S,$ and by $\{c_i\,|\,i \in I\}$ the set of exposed points of the subdifferential $\partial \varphi_2$ indexed by elements of a set $I.$ Since $\partial \varphi_1$ and $\partial \varphi_2$ are compact, the families $\{b_s\,|\,s \in S\}$ and $\{c_i\,|\,i \in I\}$ are bounded. Besides,
$$
p(x) = \sup\limits_{s \in S}\langle b_s,\,x \rangle - \sup\limits_{s \in S}\langle c_i,\,x \rangle = \sup\limits_{s \in S}\inf\limits_{i \in I}\langle b_s - c_i,\,x \rangle = \inf\limits_{i \in I}\sup\limits_{s \in S}\langle b_s - c_i,\,x \rangle\,\,\forall\,\,x \in {\mathbb R}^n.
$$
This completes the proof of the necessary part of the theorem.

The sufficiency is obvious.
\qed

Recall, that a function $p: {\mathbb R}^n \to {\mathbb R}$ is called \textit{piecewise linear} iff it is continuous and there exists a finite covering of ${\mathbb R}^n$ with solid convex cones $K_1,\,K_2,\ldots,K_s$ such that the restriction of $p$ to each cone $K_i,\,i=1,2,\ldots,s$ coincides with the restriction of some linear function to this cone.

This definition of piecewise linearity of functions goes back to Bank B. et al. \cite{Bank} (see also \cite{GorZor})

\begin{theorem}\label{th4}
A function $p:{\mathbb R}^n \to {\mathbb R}$ is piecewise linear if and only if there exists a finite two-index family $\{a_{is} \in {\mathbb R}^n\,|\,i \in I,\,s \in S\}$ which provides the following saddle representation of $p$ by linear functions:
\begin{equation}\label{e11}
p(x) = \min\limits_{i \in I}\max\limits_{s \in S}\langle a_{is},\, x \rangle = \max\limits_{s \in S}\min\limits_{i \in I}\langle a_{is},\, x \rangle, \,\,\forall\,\,x \in {\mathbb R}^n,
\end{equation}
Moreover, for each piecewise linear function $p:{\mathbb R}^n \to {\mathbb R}$ the family $\{a_{is} \in {\mathbb R}^n\,|\,i \in I,\,s \in S\}$ can be chosen in such a manner that  $\{a_{is} :=b_s - c_i \in {\mathbb R}^n\,|\,i \in I,\,s \in S\}$ where $\{b_{s} \in {\mathbb R}^n\,|\,s \in S\}$ and $\{c_{i} \in {\mathbb R}^n\,|\,i \in I\}$ are two finite families of vectors in ${\mathbb R}^n.$
\end{theorem}

{\it Proof} Due to \cite[Theorem 5.4]{GT2016} a function $p:{\mathbb R}^n \to {\mathbb R}$ is piecewise linear if and only if there exist both a finite family $\Phi$ of sublinear functions and a finite family $\Psi$ of superlinear functions which satisfy the equalities \eqref{e2} and \eqref{e3}, respectively. Repeating the arguments of the proof of Theorem \ref{th2} with the finite families $\Phi$ and $\Psi$ we construct a finite two-index family $\{a_{is},\,i \in I,\,s \in S,\}$, which provides for $p$ the saddle representation \eqref{e11}.

The converse also is true, since a function $p$ admitting the presentations \eqref{e11} is piecewise linear \cite{GorZor}.

The last statement of the theorem follows from the presentation of each piecewise linear function $p$ as a difference of two polyhedral sublinear functions \cite{Melz} (see, also, \cite[Theorem 5.3]{GT2016}),
$
p(x) = \max\limits_{s \in S}\langle b_{s},\, x \rangle - \max\limits_{i \in I}\langle c_{i},\, x \rangle, \,\,\forall\,\,x \in {\mathbb R}^n,
$
with $I$ and $S$ being finite.
\qed

\section{Conclusions}

In the paper we show that each continuous p.h. function
can be associated with a two-index family of linear functions which provides its saddle representation. This means that it admits both an inf-sup-representation and a sup-inf-representation with the same two-index family of linear functions. We also establish characteristic properties of those two-index families of linear functions which provides saddle representations of functions belonging to
the subspace of Lipschitz continuous p.h. functions as well as the subspaces of difference sublinear and piecewise linear functions.

\section*{Acknowledgements}
The research was supported by the Belarussian State Research Programm (grant ``Conversion -- 1.04.01'').


\section*{Conflict of Interest}
The authors declare that they have no conflict of interest.

\end{document}